\documentclass[11pt]{amsart}
\usepackage{amsthm,amsmath,amssymb,amscd}
\usepackage{amssymb}
\usepackage{graphics}

\newcommand{\K}{K\"{a}hler}
\newcommand{\e}{\varepsilon}
\newcommand{\del}{\partial}

\newtheorem{theorem}{Theorem}[section]

\newtheorem{proposition}{Proposition}[section]

\newtheorem{lemma}{Lemma}[section]
\newtheorem{problem}{Problem}[section]

\begin{document}

\title{On the K\"ahler classes of constant scalar curvature metrics on blow ups}

\author[Claudio Arezzo] {Claudio Arezzo}
\address{ claudio.arezzo@unipr.it \\ Universita' di Parma\\
   Parco delle Scienze 45/A\\43100 Parma\\Italy}

\author[Frank Pacard] {Frank Pacard}
\address{pacard@univ-paris12.fr\\ University Paris 12 and
Institut Universitaire de France, France}

\maketitle

\section{Introduction}

In this short paper we address the following question~:
\begin{problem}
\label{pb1}
Given a compact constant scalar curvature \K\ manifold $(M, J , g , \omega)$, of complex dimension $m := \mbox{dim}_{\mathbb C} \, M$, and having defined 
\[
\triangle : =\{(p_1, \ldots, p_n) \in M^n \quad : \quad  \exists \, a \neq b \quad p_a = p_b \},
\]
characterize the set $\mathcal{PW} = \{ (p_1,\dots , p_n, \alpha_1, \dots , \alpha_n)\}\subset (M^n \setminus \triangle) \times (0,+\infty)^n$ for which $\tilde{M} = Bl_{p_1,\dots ,p_n} M$, the blow up of $M$ at $p_1, \ldots, p_n$ has a constant scalar curvature \K\ metric (cscK  from now on) in the \K\ class
\[
\pi^*[\omega] - (\alpha_1 \, PD[E_1] +\dots + \alpha_n \, PD[E_n]),
\]
where the $PD[E_j]$ are the Poincar\'e duals of the $(2m-2)$-homology classes of the exceptional divisors of the blow up at $p_j$.
\end{problem}

This general problem is too complicated and its solution is likely to pass through the solution of  well known conjectures relating the existence of cscK  metrics with the $K$-stability of the polarized manifold.

\medskip

\noindent Yet, more specific questions are treatable and could give light also on these ambitious programs. The first natural narrowing of Problem \ref{pb1} is to require that not just one \K\ class has a cscK  representative, but that this is the case for a whole segment in the \K\ cone of $\tilde{M}$ touching the boundary at a point of the form $\pi^*[\omega]$, where $\omega$ is (necessarily) a cscK  form on $M$. Analytically this amounts to the following~:

\begin{problem}
\label{pb2}
Given a compact \K\ constant scalar curvature manifold $(M, J , g , \omega)$ characterize the set $\mathcal{APW} = \{ (p_1,\dots , p_n, a_1, \dots ,  a_n)\}\subset (M^n \setminus \triangle) \times (0,+\infty)^n$ such that $\tilde{M} = Bl_{p_1,\dots , p_n}M$ has a constant scalar curvature \K\ metric in the class 
\[
\pi^*[\omega] - \e^2 \, ( a_1 \, PD[E_1] + \dots +  a_n \, PD[E_n]) ,
\]
for all $\e$ sufficiently small. Here $\mathcal{APW}$ refers to "asymptotic points and weights", namely points and weights in this singular perturbation setting.
\end{problem}

Hence we can consider $(\alpha_1, \dots , \alpha_n)$ as an asymptotic direction in the \K\ cone for which canonical representative can be found. It is immediate to extract from \cite{P} the following ~:
\begin{theorem}
Assume that $(M, J, g, \omega)$ is a constant scalar curvature compact \K\ manifold without any nontrivial hamiltonian holomorphic vector field. Then $\mathcal{APW} = (M^n \setminus \Delta) \times (0,+\infty)^n$.
\end{theorem}

The presence of hamiltonian holomorphic vector fields greatly enhances the difficulty and the interest of the problem. In \cite{P2} the authors have attacked this problem and found an interplay between its solution and the behavior of the hamiltonian holomorphic vector fields at the $p_j$ that we briefly recall.

\medskip

First recall that  the Matsushima-Lichnerowicz Theorem asserts that the space of hamiltonian holomorphic vector fields on $(M,J, \omega)$ is also the complexification of the real vector space of holomorphic vector fields $\Xi$ which can be written as
\[
\Xi = X - i \, J \, X ,
\]
where $X$ is a Killing vector field which vanish somewhere on $M$. Let us denote by ${\mathfrak h}$, the space of  hamiltonian holomorphic vector field and by
\[
\xi_\omega : M  \longmapsto  {\mathfrak h}^*
\]
the {\em  moment} map which is defined by requiring that if $\Xi \in {\mathfrak h}$, the function $ \zeta_{\omega} : = \langle \xi_\omega  , \Xi \rangle$ is a (complex valued) Hamiltonian for the vector field $\Xi$, namely the unique solution of
\[
- \bar \del \zeta _{\omega}  = \frac{_1}{^2} \, \omega (\Xi , -) ,
\]
which is normalized by 
\[ 
\int_{M} \, \zeta_{\omega}  \, dvol_g =0 . 
\]

With these notations, the result we have obtained in \cite{P2} reads~:
\begin{theorem}
Assume that $(M, J, g, \omega)$ is a constant scalar curvature compact \K\ manifold and that $p_1, \ldots, p_n \in M$ and $a_1, \ldots, a_n >0$ are chosen so that~:

\medskip

\begin{itemize}
\item[(i)] $\xi_\omega (p_1) , \ldots, \xi_\omega (p_n)$ span ${\mathfrak h}^*$ \\[3mm]
\item[(ii)] $\sum_{j=1}^n a_j^{m-1}  \, \xi_\omega (p_j)  = 0 \in {\mathfrak h}^* $ .
\end{itemize}

\medskip

Then, there exist $\e_0 > 0$ such that, for all $\e \in (0, \e_0)$, there exists on $\tilde M = Bl_{p_1,\dots ,p_n}M$ , a constant scalar curvature \K\ metric $g_\e$ associated to the \K\ form
\[
\omega_\e \in \pi^{*} \, [\omega] - \e ^{2} \, (a_{1, \e} \, PD[E_1]
+ \ldots + a_{n, \e} \, PD[E_n]),
\]
where
\begin{equation}
|a_{j ,\e} - a_j | \leq c \, \e^{\frac{2}{2m+1}} .
\label{disto}
\end{equation}
Finally, the sequence of metrics $(g_\e)_\e$ converges to $g$ in ${\mathcal C}^{\infty} (M \setminus \{p_1, \ldots, p_n\})$. \label{th:2}
\end{theorem}
Therefore, in the presence of nontrivial hamiltonian holomorphic vector fields, the number of points which can be blown up, their position, as well as the possible \K\ classes on the blown up manifold have to satisfy some constraints.

\medskip

It is not hard to see from the proof in \cite{P2} that the Mapping
\[
(a_1, \ldots, a_n) \longmapsto (a_{1, \e}, \ldots, a_{n, \e})
\]
is continuous. Indeed, this follows from the construction itself which only uses fixed point theorems for contraction mappings and hence the metric we obtain depends smoothly on the parameters of the construction.

\medskip

Theorem~\ref{th:2} has two major drawbacks : First, we lose control on the \K\ classes on $\tilde M$ for which constant scalar curvature \K\ metrics can be constructed, second there are severe restrictions on the set of points and asymptotic directions. 

\medskip

The key idea to fill these gaps is to note that the construction of \cite{P2} is in fact a construction of the Riemannian metric $g_\e$ and this is reflected by the fact that the sequence of metrics constructed converges to the initial metric $g$ and also in the fact that condition (ii) really depends on the choice of the metric $g$.  

\medskip

Now, on the one hand, the origin of (ii) stems  from the existence of hamiltonian holomorphic vector fields on $(M, J)$ and in fact (ii) imposes on the choice of the asymptotic directions $(a_1, \ldots, a_n)$ as many constraints as the dimension of ${\mathfrak h}$. 

\medskip

On the other hand, the existence of hamiltonian holomorphic vector fields is also related to the non-uniqueness of the constant scalar curvature \K\ metric on $M$.  More precisely, ${\mathfrak h}$ is the Lie algebra of the group of automorphisms of $(M, J, g, \omega)$ and as such also parameterizes near $g$ the space of constant scalar curvature \K\ metrics in a given \K\ class $[\omega]$ and for a given scalar curvature. Observe that this space has dimension  equal to $\mbox{dim} \, {\mathfrak h}$. Therefore, we can Apply the result of Theorem~\ref{th:2} not only to the metric $g$ itself but also to the pull back of $g$ by any biholomorphic transformation. 

\medskip

Since condition (ii) depends on the choice of the metric, if we are only interested in the \K\ classes on the blown up manifold, we get more flexibility in the choice of the asymptotic parameters (observe that the dimension of the space of constant scalar curvature \K\ metrics near $g$ (with fixed scalar curvature) is precisely equal to the number of constraints on the choice of the asymptotic parameters). This observation allows us to complement the result of Theorem~\ref{th:2} and get the~:
\begin{theorem} Assume that $(M, J, g, \omega)$ is a constant scalar curvature compact \K\ manifold and that $p_1, \ldots, p_n \in M$ and $a_1, \ldots, a_n >0$ are chosen so that~:

\medskip

\begin{itemize}
\item[(i)] $\xi_\omega (p_1) , \ldots, \xi_\omega (p_n)$ span ${\mathfrak h}^*
$ \hspace{4cm} (genericity condition)\\[3mm]
\item[(ii)] $\sum_{j=1}^n a_j^{m-1} \, \xi_\omega (p_j)  = 0 \in {\mathfrak h}^* $
\hspace{4cm} (balancing condition)
\\[3mm]
\item[(iii)] no element of ${\mathfrak h}$ vanishes at every point $p_1, \ldots , p_n$. \hspace{1cm} (general position condition)
\end{itemize}

\medskip

Then $(p_1, \dots, p_n , a_1, \dots, a_n) \in\mathcal{APW}$.
\label{th:3}
\end{theorem}
Therefore, we can indeed prescribe the exact value of the asymptotic direction in which the \K\ classes in perturbed at the expense of  imposing that no hamiltonian holomorphic vector field vanishes at  every point we blow up. 

\medskip

The genericity condition  is purely technical and it does not seem to hide any deep geometric nature. Indeed, as observed in \cite{P2}~:
\begin{lemma}
With the above notations, assume that $n \geq \mbox{dim} \, {\mathfrak h}$. Then, the set of points $(p_1, \ldots, p_n) \in M^n \setminus \triangle$ satisfying the genericity condition is open and dense. \label{le:1}
\end{lemma}

The balancing condition is certainly the heart of the problem, encoding the relevant stability property of $\tilde{M}$. For example when all the $a_j$ are rationals, the balancing condition is easily translated in the Chow polystability of the cycle $\sum_j a_j^{m-1} \, p_j$ with respect to the action of the automorphism group of $M$.

\medskip

In a remarkable recent paper Stoppa \cite{st} has proved, among other things, that it if the cycle  $\sum_j a_j^{m-1} \,  p_j$ is Chow unstable, then  $(p_1, \dots, p_n, a_1, \dots, a_n )$ does not lie in $\mathcal{APW}$. With a beautifully careful algebraic analysis he has in fact related  a destabilizing configuration for the points to a destabilizing  configuration of the blown up manifold giving a quantitative measure of the reciprocal unstabilities.

\medskip

Going back to our problem, we first observe that the combination of the three above condition still leaves flexibility in the choices~:
\begin{theorem}
With the above notations, assume that $n \geq \mbox{dim} \, {\mathfrak h}+1$ then, the set of points $((p_1, \ldots, p_n), (a_1, \ldots, a_n))   \in (M^n \setminus \triangle)  \times (0, \infty)^n $ such that condition (i), (ii) and (iii) are fulfilled is open in $(M^n \setminus \triangle) \times (0, \infty)^n$. 
\label{le:2}
\end{theorem}
Openness in the choice of the points was already contained in \cite{P2}. What we will prove in this short pPer is the openness in the choice of the asymptotic directions.

\medskip

We can better understand the topology of $\mathcal{APW}$ by looking at 
\[
\begin{array}{rcccr}
&  & \mathcal{APW} &  &  \\[3mm]
& \pi_1 \swarrow  &   &  \searrow  \pi_2&\\[3mm]
M^n \setminus \triangle&  & & &  (0, \infty)^n  
\end{array}
\]
 and define 
 \[
 \mathcal{AP}=\pi_1(\mathcal{APW})Ê\qquad \mbox{and}  \qquad \mathcal{AW}=\pi_2(\mathcal{APW}).
 \]
 
With these notations, we obtain~:   
\begin{theorem} 
 \label{thP}
Assume that $(p_1, \dots, p_n) \in \mathcal{AP}$ and further assume that the general position condition holds, then
 \newline
 $\pi_2(\pi_1^{-1}(p_1, \dots, p_n))$ is an open (nonempty) subset of $(0, \infty)^n$.
 \end{theorem}
 
 And we also have the~:
 \begin{theorem}
 \label{thaw}
Assume that $(a_1, \dots, a_n) \in \mathcal{AW}$ and further assume that there exists $(p_1, \dots, p_n, a_1, \dots, a_n)    \in \pi_2^{-1}(a_1,\dots, a_n)$ for which the general position condition holds, then  $\pi_1(\pi_2^{-1}(a_1, \dots, a_n))$ is an open dense subset of $ M^n \setminus \triangle$.
 \end{theorem}
Hence the general position condition shows that, by moving the cscK  representative in $[\omega]$, the balancing condition is a very flexible one.

\medskip
 
Theorems \ref{thP} and \ref{thaw} are of completely different nature. Theorem \ref{thP}  follows from an implicit function argument Applied to the set of solutions of the balancing condition, and it is of a local nature. On the other hand Theorem \ref{thaw} follows from a  suitable interpretation of the balancing condition in terms of the geometry of moment mPs. In this language, we can interpret $(p_1, \dots, p_n)\in \pi_2^{-1} (a_1, \dots, a_n)$ as a point in the zero set of a natural moment mP, and the general position condition is readily translated in the fact that this point is regular. The general theory then provides openness and density of the orbits of $(p_1, \dots, p_n)$ through the action of the automorphisms group, which in turns implies the result.

\medskip
 
We should note in this regard that it is a hopeless and confusing task to check conditions (i), (ii) and (iii) when the points move in these orbits but one should simply transport the solution associated to $(p_1, \dots, p_n, a_1, \dots, a_n)$ on $\tilde{M}$ to a solution at  $( \mathfrak g (p_1), \dots,  \mathfrak g (p_n), a_1, \dots, a_n)$ on $ \mathfrak g (\tilde{M})$ where $ \mathfrak g \in Aut(M)$.

\medskip
 
It is important to emphasize that Theorem \ref{thP} cannot be improved to get a density result, as 
Stoppa \cite{st} has found explicit bounds for the choice of weights to have cscK   metrics in the blow ups of even deceptively simple examples as the projective plane  (this result has then been strengthened by Della Vedova \cite{dv} to encompass the case of extremal metrics).

\medskip

We conclude this pPer by analyzing $\mathcal{APW}$ in the special case of manifolds for which  $dim \, {\mathfrak h} =1$ and $n=2$. We show that in this case  $\pi_2(\mathcal{APW})=(0, \infty)^2$. Recall that, among these type of manifolds, there are nontrivial explicit cscK  metrics obtained by LeBrun \cite{lb} by the so-called {\em moment construction}.  In this case we can also characterize $\pi_1(\mathcal{APW})$. 

\medskip

The analysis carried through in this note can be adPted to analyze the similar problem for extremal metrics in the sense of Calabi. In this case the algebraic analysis done by Stoppa \cite{st} for $K$-stability has been completed by Della Vedova \cite{dv} for the {\em relative} $K$-stability in the sense of Szekelyhidi \cite{sz} .

\section{Proof of the results}

We now proceed with the proof of the different results. 

\medskip

\subsection{Proof of  Theorems~\ref{th:3} and \ref{thP}.} Given a holomorphic vector field $\Xi \in {\mathfrak h}$  (close to $0$) we consider $\phi_\Xi$ to be the bi-holomorphic transformation obtained by exponentiating the vector field $\Xi$, namely we
consider the flow of the vector field $\Xi/||\Xi||$ at time $t= ||\Xi||$. Here $||\, \cdot\, ||$ is any norm on ${\mathfrak h}$, they are all equivalent since this space is finite dimension. Observe that we have
\[
[\phi_{\Xi}^* \, \omega] = [\omega],
\]
and 
\[ 
{\bf s} (\phi^*_\Xi  \, \omega ) = {\bf s}(\omega) .
\]
Therefore,  for all $\Xi \in \mathfrak h$, the \K\ form  $\phi_\Xi^* \, \omega$ can be used in Theorem~\ref{th:2} to construct constant scalar curvature \K\ metrics on the blow up of $M$ at $p_1, \ldots, p_n$ in a \K\ class close to
\[
\pi^{*} \, [\omega] - \e ^{2} \, (a_{1} \, PD[E_1] + \ldots + a_{n} \, PD[E_n]).
\]

The above discussion shows that we should be interested in the set \K\ forms $\tilde \omega \in [\omega]$, with ${\bf s} (\tilde \omega ) = {\bf s} (\omega)$, points $(p_1, \ldots, p_n) \in M^n \setminus \triangle$ and asymptotic directions $(a_1, \ldots, a_n) \in (0, \infty)^n$ solution of the equation
\[
a_1^{m-1} \, \xi_{\tilde \omega} (p_1) + \ldots + a_n^{m-1} \, \xi_{\tilde \omega} (p_n) = 0 \in
{\mathfrak h}^*.
\]
Let us assume that we have such a solution $(a_1, \ldots, a_n) \in (0, \infty)^n$ and $(p_1, \ldots , p_n) \in M^n \setminus \triangle$ for the \K\ form $\tilde \omega = \omega$ itself. We would like to know if, close to this solution, we can move freely the coefficients $a_j$ and the points $p_j$. To this aim, we consider the Mapping
\[
{\mathfrak S}  : (M^n \setminus \triangle ) \times (0, +\infty)^n  \times {\mathfrak h} \longrightarrow {\mathfrak h}^*
\]
defined by
\[
{\mathfrak S}  ((b_1, \ldots, b_n), (q_1, \ldots, q_n) , \Xi) = b_1^{m-1}  \, \xi_{\tilde \omega} (q_1) + \ldots + b_n^{m-1}  \, \xi_{\tilde \omega} (q_n),
\]
where 
\[ 
\tilde \omega  =  \phi^*_\Xi \, \omega .
\]
Given $[\omega]$, the possible points which can be blown up and the possible asymptotic directions in which the \K\ class $[\omega]$ can be perturbed are just the projection over the first two components of the zeros of the Mapping ${\mathfrak S}$.

\medskip

Let us assume that the differential of ${\mathfrak S}$ with respect to $\Xi$, computed at the point $((p_1, \ldots, p_n), (a_1, \ldots , a_n), 0) \in (M^n \setminus \triangle)  \times (0, \infty)^n \times \mathfrak h$, is an isomorphism between ${\mathfrak h}$ and ${\mathfrak h}^*$. The implicit function theorem Applied to the Mapping ${\mathfrak S}$ guaranties that, close to $((p_1, \ldots, p_n) , (a_1, \ldots, a_n),  0)$, the set of solutions of
\begin{equation}
{\mathfrak S}  (  (q_1, \ldots, q_n) , (b_1, \ldots, b_n), \Xi)  = 0 , \label{eq:1}
\end{equation}
is parameterized by $ (q_1, \ldots, q_n)$ and $(b_1, \ldots, b_n)$. In other words,  equation (\ref{eq:1}) can be solved as
\[
\Xi = HV  ((b_1, \ldots, b_n), (q_1, \ldots, q_n)) ,
\]
for some (smooth) Mapping $ HV $ defined from a neighborhood of $( (p_1, \ldots , p_n), (a_1, \ldots, a_n))$ in $(M^n \setminus \triangle ) \times (0, \infty)^n$ into $\mathfrak h$ and satisfying $HV ((p_1, \ldots, p_n), (a_1, \ldots, a_n)) = 0$. 

\medskip

To complete the proof of Theorems~\ref{th:3} and \ref{thP}, we consider the Mapping
\[
(b_1, \ldots, b_n) \longmapsto HV  ((p_1, \ldots, p_n), (b_1, \ldots, b_n))
\]
which is defined in a neighborhood $U$ of $(a_{1}, \ldots, a_n)$, with values in ${\mathfrak h}$. Observe that, by construction
\[
\sum_{j=1}^n \, b_j^{m-1} \, \xi_{\phi_\Xi^* \, \omega} (p_j) = 0
\]
if $\Xi = HV ((p_1, \ldots, p_n), (b_1, \ldots, b_n))$. Moreover, reducing $U$ if this is necessary, we can always assume that
\[
\xi_{\phi_\Xi^* \, \omega} (p_1) , \ldots, \xi_{\phi_\Xi^* \,  \omega} (p_n)
\]
span ${\mathfrak h}^*$, since this is true when $\Xi =0$ thanks to condition (i). Hence, we can use the result of Theorem~\ref{th:2} aApplied to the metric associated to $\phi^*_{\Xi} \, \omega$, the points $p_1, \ldots, p_n$ and the coefficients $(b_1, \ldots, b_n) \in U$. This provides a \K\ metric in the \K\ class
\[
\pi^{*} \, [\omega] - \e ^{2} \, (b_{1, \e}^{m-1} \, PD[E_1] + \ldots + b_{n, \e}^{m-1} \, PD[E_n]),
\]
where the coefficients $b_{j, \e}$ depend (smoothly) on the points $b_1, \ldots, b_n$. To summarize, we have defined  a Mapping
\[
C_\e :  (b_1, \ldots , b_n) \in U \longmapsto  (b_{1, \e} , \ldots, b_{n, \e}) \in (0,\infty)^n .
\]
As already mentioned, this Mapping is at least continuous and is close to the identity since (\ref{disto}) implies that
\[
|| C_\e ( (b_1, \ldots , b_n)) -(b_{1} , \ldots, b_{n}) || \leq c \, \e^{\frac{2}{m+2}}
\]
Clearly, $\mbox{deg} (C_\e , (a_1, \ldots, a_n) ; U) =1$ for $\e$ small enough and hence the image of $U$ by $C_\e $ contains an open neighborhood of $(a_1, \ldots, a_n)$, provided $\e$ is chosen small enough. This implies that, for all $\e$ small enough, it is possible to find $(b_1, \ldots, b_n) \in U$ so that
\[
C_\e ( (b_1, \ldots , b_n))  = (a_1 , \ldots, a_n) ,
\]
and this completes the proof of Theorems~\ref{th:3} and \ref{thP}.

\medskip

Therefore, the only thing which remains to be understood is when the differential of ${\mathfrak S}$ with respect of $\Xi$ is an isomorphism.

\medskip

\noindent {\bf The differential of ${\mathfrak S}$ with respect of $\Xi$} Let $\Xi \in {\mathfrak h}$ be given and $t \in {\mathbb R}$. We consider the one dimensional family of \K\ forms
\[
\omega_t = \phi^*_{t \Xi} \, \omega.
\]
First observe that, for $t$ small we can expand
\[
\omega_t = \omega + i \partial \, \bar \partial \, (t \, f) + {\mathcal O} (t^2) ,
\]
where $f$ is precisely the (real valued) potential associated to $\Xi$ given by
\[
- \bar \partial \, f  =  \tfrac{1}{2} \, \omega (\Xi , -) .
\]
Recall that we can write
\[
\Xi = X - i \, J \, X ,
\]
for some Killing vector field $X$ (for the metric $g$ associated to $\omega$) and we can also write
\begin{equation}
- d \, f  =  \omega (X , -) . 
\label{eq:4}
\end{equation}

Observe that if we consider the metric $\omega_t$,  any fixed holomorphic vector field $\tilde \Xi  \in {\mathfrak h}$ is associated to a (complex valued) potential (depending on $t$), which is defined by
\begin{equation}
- \bar \partial  \tilde f_t  = \tfrac{1}{2}Ê\, \omega_t ( \tilde \Xi , -)
\label{eq:zs}
\end{equation}
and  $\langle \tilde \xi_{\omega_t} , \tilde \Xi \rangle = \tilde f_t$, where $\xi_{\omega_t}$ is the moment mP associated to the \K\ form $\omega_t$.

\medskip

Differentiating (\ref{eq:zs}) with respect to $t$, at $t=0$, we find
\begin{equation}
- \bar \del \langle \dot \xi_{\omega}, \tilde \Xi \rangle =  \tfrac{i}{2} \, \del \bar \del f ( \tilde \Xi , -)
\label{eq:qs}
\end{equation}
where  $\dot \xi_\omega$ is the first variation of $f \longmapsto \xi_{\omega + i \del \bar \del f}$.  Working in local coordinates and using the fact that $\Xi$ is holomorphic, one checks that the right hand side of (\ref{eq:qs}) is equal to $\frac{i}{2}Êd \bar z^a \frac{\partial}{\partial \bar z^a} \left( \tilde  \Xi^b  \, \frac{\del f}{\del z^b}\right)$. Hence, we conclude that 
\[
\langle \dot \xi_{\omega}, \tilde \Xi \rangle =  - \tfrac{i}{2} \,\tilde \Xi \, f .
\]

It is enough to consider the set of holomorphic vector fields $\Xi_2$ which can be written as
\[
\tilde \Xi = \tilde X - i \, J \, \tilde X ,
\]
for some Killing vector field $\tilde X$ (for the metric $g$). Given the definition of $f$, we get
\begin{equation}
\begin{array}{rllll}
\langle \dot \xi , \tilde \Xi \rangle & = & - \displaystyle \frac{_i}{^2} \, \, \tilde \Xi \, f \\[3mm]
&  = & - \displaystyle  \frac{_i}{^2} \, df (\tilde \Xi) \\[3mm]
& = & - \displaystyle  \frac{_i}{^2}  \, df ( \tilde X - i \, J \, \tilde X) \\[3mm]
& = & \displaystyle  \frac{_i}{^2}  \, \omega ( X , \tilde X - i \, J \, \tilde X) \\[3mm]
& = & \displaystyle  \frac{_1}{^2}  \, \left( g (   X ,  \tilde X) + i \, g( J\, X,  \tilde X) \right)
\label{diffmom}
\end{array}
\end{equation}
The important point is that
\[
(\tilde \Xi ,  \Xi )_{Her} : = g ( X , \tilde X)  + i \, g \, (J \, X,  \tilde X)
\]
is a positive definite Hermitian form. Alternatively, this corresponds to
\[
( \tilde \Xi , \Xi  )_{Her} : =   - \frac{_i}{^2} \,  \omega ( \tilde \Xi , \bar \Xi ).
\]

We denote by $L$ the differential of ${\mathfrak S}$ with respect to $\Xi$, computed at $(a_1, \ldots , a_n)$, $(p_1, \ldots, p_n)$ and $\Xi=0$. So that
\[
L : {\mathfrak h}\longrightarrow {\mathfrak h}^*.
\]
and $L( \Xi) \in {\mathfrak h}^*$. Using the above computation, we conclude that
\begin{equation}
L(\Xi) = \tfrac{1}{2}Ê\, \sum_{j=1}^n \, a_j^{m-1}  \, ( - , \Xi)_{Her} (p_j) .
\label{eq:3}
\end{equation}
Now,  $L$ generates a positive Hermitian form on $\mathfrak h$ by
\[
(\tilde X, \Xi ) = \tfrac{1}{2} \, \sum_{j=1}^n \, a_j^{m-1}  \, (\tilde , \Xi)_{Her} (p_j) .
\]
Clearly, this form is non degenerate if and only if there is no holomorphic vector field $\Xi \in \mathfrak h$ which vanishes at every  point $p_1, \ldots, p_n$ (therefore, we need $n \geq \mbox{dim}Ê\mathfrak h$). To summarize, we have proved the~: 
\begin{proposition} Assume that there are no nontrivial element of ${\mathfrak h}$ vanishing at  every $p_1, \ldots, p_n$, then the differential of ${\mathfrak S}$ with respect to $\Xi$, computed at $((a_1, \ldots, a_n), (p_1, \ldots, p_n),0)$ is an isomorphism from ${\mathfrak h}$ into ${\mathfrak h}*$.
\end{proposition}

All the above discussion seems to point out that the really important object is the zero set of the mapping ${\mathfrak S}$, or more precisely, its projection over the first two entries (the set of points which can be blown up and the set of asymptotic directions toward which the \K\ class can be deformed). This also explains the role of the zeros of the holomorphic vector field, role which was completely occulted in \cite{P2} since the only important condition was associated to the potentials not the gradient of the potentials.

\subsection{ Proof of  Theorem~\ref{thaw}.} We consider the action of the hamiltonian isometry group $H$ for some metric $g$. We can also consider $H$ acting on $M^n$ equipped with the weighted metric
\[
a_1^{m-1} \, g + \ldots + a_n^{m-1} \, g .
\]
The moment mP for this action is then given by
\[
\mu : = a_1^{m-1} \, \xi_\omega + \ldots + a_n^{m-1} \, \xi_\omega .
\]
In our case, Theorem 7.4 in \cite{ki} asserts that, if $\mu^{-1} (0) \neq \emptyset$  and if there exists $(p_1, \dots, p_n) \in \mu^{-1}(0)$ satisfying the general position condition, then $(H \otimes {\mathbb C}) \, \cdot \mu^{-1} (0) $ is open and dense in $M^n \setminus \triangle$.

\medskip

In other words, if we have a set of points $p_1, \ldots, p_n$ for which 
\[ 
\sum_{j=1}^n \, a_j^{m-1} \, \xi_\omega (p_j) = 0,
\]
so that $(p_1, \ldots, p_n ) \in \mu^{-1} (0)$ then the action of $H \otimes {\mathbb C}$, the complexification of $H$, provides an open dense set of points ${\mathcal U} \subset M^n \setminus \triangle$ for which the condition
\[
\sum_{j=1}^n \, a_j^{m-1} \, \xi_{\phi^*\omega} (q_j) =0 ,
\]
is fulfilled for some automorphism $\phi$ (depending on $(q_1, \ldots, q_n )\in {\mathcal U}$).  Now $\phi$ lifts to a biholomorphic mP $\tilde{\phi}\colon Bl_{p_1,\dots,p_n}M \rightarrow Bl_{q_1,\dots,q_n}M$, and since we know that there exists a family of cscK  forms $\omega_{\e} \in \pi^*[\omega] - \e ^2 \, (a_1 \, PD[E_1]+\dots + a_n \, PD[E_n])$ on  $Bl_{p_1,\dots,p_n} \, M$,  the family  $(\tilde{\phi}^{-1})^{*}(\omega_{\e})$ is the seeked family of cscK  forms on $Bl_{q_1,\dots,q_n}M$. This completes the proof of the result.

\section{An important example}

Let us consider the simplest case where ${\mathfrak h} = \mbox{Span} \{\Xi\}$ and where we want to blow up $2$ points. Then the condition on the points and the asymptotic directions which can be considered for the blow up procedure, becomes 
\[
a_1^{m-1} \, \xi_\omega (p_1) + a_2^{m-1} \, \xi_\omega (p_2) =0
\]
If
\[
\zeta : = \langle \xi_\omega, \Xi\rangle
\]
this is just
\begin{equation}
a_1^{m-1} \zeta (p_1) + a_2^{m-1} \, \zeta (p_2) =0 .
\label{ball}
\end{equation}
Therefore, according to the result of \cite{P2} we can blow up any two points $p_1, p_2$ with directions $a_1, a_2$ satisfying (\ref{ball}), provided $\zeta (p_1)$ and $\zeta (p_2)$ are not zero and have different signs. At this point, it looks like there is a constraint on the directions ! However, as a consequence of our result, we see that we can locally move the coefficients $a_1, a_2$ freely provided the vector field $\Xi$ does not vanish at the two points $p_1, p_2$. Indeed, the formula that guarantees that $L$ is invertible just reduces to
\[
L(\Xi, \Xi) \neq 0 
\] 
which in this simple case reads
\[
a_1^{m-1} \, (\Xi, \Xi)_{Her} (p_1) + a_2^{m-1} \, (\Xi, \Xi)_{Her} (p_2) \neq 0 .
\]
Therefore, if $\Xi$ does not vanish at both $p_1$ and $p_2$ (i.e. the general position condition is satisfied), then we can use the blow up procedure of \cite{P2} and we see that the set of directions of deformation of the K\"ahler classes is open. Observe that we obtain uniqueness of the corresponding constant scalar curvature K\"ahler metric since blowing up the points has "killed" the unique holomorphic vector field on $M$.

\medskip

Observe that given $a_1, a_2 >0$ it is always possible to find points $p_1, p_2 \in M$ satisfying (\ref{ball}), i.e.  $\pi_2(\mathcal{APW})=(0,+\infty)^2$. Indeed, denote
\[
a^- = \min \zeta < 0 < \max \zeta = a^+ .
\]
(Recall that the average of $\zeta$ over $M$ is zero). The intermediate value Theorem shows that one can find points $p_1, p_2$ such that
\[
a^- < \zeta(p_1) =  \frac{a^- a^{+} \, a_2^{m-1}}{\sqrt{(a^- a_1^{m-1})^2 + (a^{+} a_2^{m-1})^2}}   < 0 ,
\]
and 
 \[
 0< \zeta(p_2) = \frac{-  \, a^+ \, a^{-} a_1^{m-1}}{\sqrt{(a^- a_1^{m-1})^2 + (a^{+} a_2^{m-1})^2}} < a^+,
\]
and hence (\ref{ball}) is satisfied for these two points.

\medskip

The situation just described is far from an artificial speculation. Among these manifolds having only one holomorphic vector field vanishing somewhere, there is a well known class of examples of \K\ constant scalar curvature manifolds which are neither  products nor Einstein, and have been discovered and investigated by Lebrun in \cite{lb} and also in \cite{lbs}, \cite{klbp}. 

\medskip

Let us recall that such surfaces are blow ups at finite set of points along the zero section of manifolds of the type $P (\mathcal L \oplus\mathcal{O})$, where $\mathcal L$ is a line bundle  of positive degree over a Riemann surface of genus greater than $1$. Such a procedure makes only the Euler vector field $\Xi$ survive on $M$.  

\medskip

In this last setting, we define $\phi_t$ to be the flow associated to the Euler vector field $\Xi$. Given a point $p$ we write $p(t) = \phi_t (p)$ to be the image by flow of $\Xi$ passing through $p$ at time $t$. We have 
\[
\frac{d}{dt} \zeta {(p(t))} =  d\zeta (\Im \Xi) =  - d_\zeta  (J \, X) = \omega (X, JX) = g(X,X),
\]
by definition of $\zeta$. Observe that the flow $\phi$ preserves the fibers and thanks to this formula, the function $t \longmapsto \zeta (p(t))$ is monotone increasing with $t$. 

\medskip

Now, if $p_1, p_2$ belong to the same fiber and do not belong to the zero section nor to the infinity section,  given $a_1, a_2 >0$ one can find $t \in \mathbb R$ such that 
\[
a_1\, \zeta (p_1(t)) + a_2 \, \zeta (p_2(t)) = 0 .
\]
Hence, according to Theorem~\ref{th:3}, we can blow up $(M, J, g, \omega)$ at the points $p_1(t)$ and $p_2(t)$ and find a cscK  metric in the \K\ class corresponding to the weights $a_1$ and $a_2$.  This is clearly equivalent to produce on the blow up of $(M, J, \phi^*_t g, \phi^*_t \omega)$  at the points $p_1, p_2$, a cscK  metric in \K\ class corresponding to the weights $a_1$ and $a_2$. 

\medskip

The metrics we consider have the remarkable property that $\omega (\Xi, \bar \Xi)$, on each fiber of the line bundle, does not depend on the point chosen on the level set of the function $\zeta$.  Since the image of a level set of $\zeta$ by the flow $\phi$ is another level of  $\zeta$ and each fiber is preserved by the flow, it is enough to choose the points $p_1$ and $p_2$ at different "moment heights" for the above discussion to hold. 
Since the genericity condition is  obviously satisfied, this shows that $\pi_1(\mathcal{APW})$ contains $M^2\setminus \mathcal{M}$, where 
\[
\mathcal{M} = \{ (p_1,p_2) \quad : \quad  \zeta_{\omega}(p_1)=\zeta_{\omega}(p_2), \, \mbox{for some (hence {\em any}) cscK  metric}\},
\] 
which is clearly open and dense in $M^2 \setminus \triangle$.

\medskip

Conversely, if $(p_1,p_2) \in \mathcal{M}$, but not on the zero or infinity section of $\mathcal L \oplus \mathcal{O}$, then the balancing condition and the genericity condition cannot be simultaneously satisfied for any cscK  metric on $M$, hence we cannot conclude that $(p_1,p_2)$ lies in $\mathcal{AP}$. 

\medskip

The last case to analyze is when both $p_1, p_2$ both lie on the zero section or on the infinity one. In this case, the idea is to work equivariantly as in \cite{Ps} with respect to $K$, the group of isometries generated by $\Re \, \Xi$, and obtain extremal metric on $\tilde M$. This time we obtain extremal K\"ahler metrics without any constraint on the possible asymptotic directions toward which the \K\ class $\pi^* [\omega]$ can be deformed. But by Proposition 2.1 in \cite{Ps} these metrics cannot be cscK  (since the balancing conition is not satisfied) and hence these \K\ classes do not have cscK  representatives by \cite{ca2}.

\medskip

{\bf Acknowledgment  :} The authors would like to thank the referee for useful comments on the last section of the paper.  The second authors would like to acknowledge the hospitality and support of the Forschungsinstitut f\"ur Mathematik, at ETH Z\"urich, where this paper was written.

\end{document}